\documentclass[11pt,letterpaper,oneside]{article}
\usepackage{amsmath,amsthm,amssymb,amscd,abbrViro,graphicx,xcolor}
\graphicspath{  {./figs/} }
\usepackage{pgfplots,tikz-cd,tikz}
\begin{document} 
\title{Fundamental group in the projective knot theory}
\author{Julia Viro and Oleg Viro}

\date{}
\maketitle
\vspace{1cm}

\hspace*{5.2cm}{To Volodia Turaev with gratitude}\\ 
\rightline{for lifelong friendship and inspiration}
\vspace{1cm}

\begin{abstract} \noindent In this paper, properties 
of a link $L$ in the projective space $\rpt$ are  related to properties of
its group 
$\pi_1(\rpt\sminus L)$:\\ 
\hspace*{8pt}{\it $L$ is isotopic to a projective line if and only if $\pi_1(\rpt\sminus
L)=\Z$.} \\
\hspace*{8pt}{\it $L$ is isotopic to an affine circle if and only if $\pi_1(\rpt\sminus
L)=\Z*\ZZ$.} \\
\hspace*{8pt}{\it $L$ is isotopic to a link disjoint from a projective
plane if and only if $\pi_1(\rpt\sminus L)$ contains a non-trivial element of order
two.}\\
A simple algorithm which finds a system of generators and relations for 
$\pi_1(\rpt\sminus L)$ in terms of a link diagram of $L$ is provided. 
\end{abstract}

\section{Introduction}\label{s1}

\subsection*{Projective knot theory}\label{s1.1}
A {\em classical link\/} is a smooth closed 1-submanifold of $\R^3$.
If a link is connected, then it is called a {\em knot.\/}  
The {\em classical knot theory\/} studies classical knots and links.
Many notions and results of the classical knot theory extend to knots 
and links in other 3-manifolds. 
Any closed connected 3-manifold can be presented as 
a compactification of {$\R^3$}. The 3-sphere $S^3$ is the result of
adding a single point to $\R^3$. This addition is inessential for topology
of knots and links. Knots and links in $S^3$ are objects of the same
classical knot theory.  

A substantial enlargement for the set of knots and links happens when we
pass to the next simplest ambient 3-manifold, which is 
the real projective space $\rpt$.
It can be obtained from $\R^3$ by adding ``points at infinity'' as any 
other closed connected 3-manifold, but $\rpt$ is the only closed 3-manifold 
for which the added points are all alike and constitute a surface without 
singularities. The projective space has many other special features, which
provide opportunities which are not available for studying links in
other 3-manifolds. 

On the other hand, the projective space appears naturally outside of 
topology, and this gives extra reasons for study of links in the projective 
space.
 
In the real algebraic geometry the projective space is even more profound 
than the Euclidean 3-space or 3-sphere.  Recently there was a significant 
progress in real algebraic knot theory, which studies knots that are real 
algebraic curves. See \cite{MO} and \cite{Bjork}.
In that area, projective space is the most elementary and natural ambient 
space.

Any straight line in $\R^3$ gives rise to a projective line in $\rpt$. 
Topologically, a projective line is a circle. 
A configuration of lines in $\R^3$ skew to each other gives rise to a 
link in $\rpt$. See \cite{DV}. 
 
A {\em projective link} $L$ is a smooth closed 1-submanifold  of $\rpt$. 
If $L$ is connected, then it is called a {\em projective 
knot\/}. Projective knots and links are studied in 
the {\em projective knot theory\/}. It occupies a special place between
the classical knot theory and theory of links in arbitrary 3-manifolds.

This paper  appears as a fragment of a book by both authors on the projective 
link theory. It is an expanded version of a preprint \cite{V} of the second author.

Some parts of the classical knot theory have no clear generalizations
to knots in arbitrary ambient 3-manifold, but have well-developed 
counter-parts in the projective knot theory. 

For example, a large part of the classical knot theory deals with 
{\em link diagrams\/}, decorated generic projections of a link to a plane
that are used to describe classical links up to a smooth isotopy. 

\begin{figure}[h]
\centerline{
\includegraphics[scale=.65]{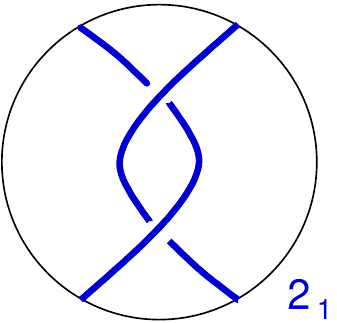}\hspace*{1cm}
\includegraphics[scale=.65]{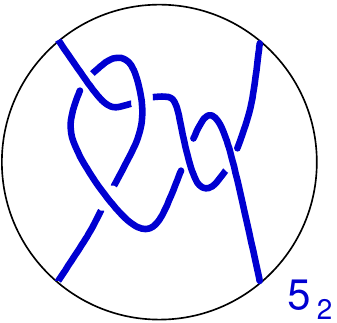}\hspace*{1cm}
\includegraphics[scale=.65]{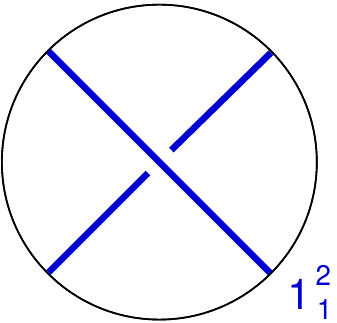}\hspace*{1cm}
\includegraphics[scale=.65]{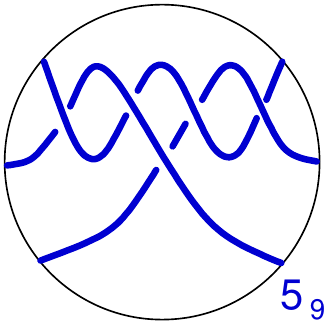}}
\vspace{10pt}
\caption{Projective link diagrams}
\label{f1}
\end{figure}

Projective links also can be described up to smooth isotopy by pictures 
similar to link diagrams that are used to describe classical links.
These pictures are called {\em projective link diagrams.\/} 
Projective link diagrams differ from classical link diagrams
in the way that a projective link diagram is placed not in the
plane, but in a disk; arcs of the diagram meet  the boundary of the
disk in pairs of antipodal points. Examples of projective link diagrams 
are shown in Figure \ref{f1}. For more details see \cite{D1}.

Most invariants of classical links can be easily calculated in terms of 
link diagrams. For some invariants, like the Jones polynomial and 
Kauffman bracket polynomial, the very definitions rely
on diagrams. The first author used projective link diagrams for generalizing
these classical links invariants to projective links, see \cite{D1}.

\subsection*{Fundamental groups in the knot theory}\label{s1.2} 
The fundamental group plays a central role in topology of 3-manifolds.
One of the most powerful topological invariants of a classical link 
$L$ is the fundamental group $\pi_1(\R^3\sminus L)$ of the complement. 
If $L$ is a knot, then $\pi_1(\R^3\sminus L)$ is even 
called just the {\em knot group.\/} 
(The term {\it link group\/} is used for a certain quotient group 
of $\pi_1(\R^3\sminus L)$, see \cite{Milnor}.)

Many geometric problems of knot theory have been solved in terms of 
knot groups. For example, in 1910 Dehn \cite{DH} proved that 
a knot is isotopic to the unknot if and only if its group is isomorphic to
$\Z$ (there was a gap in the proof, which was repaired by Papakyriacopoulus
in 1956). 

The power of the knot group  is compromised by algorithmic
difficulties in operating with non-commutative groups. To the best of our
knowledge, there is no algorithmic problem in the classical knot theory,
which have been solved via reducing  to a problem about the knot group (with
possibly extra structures like group system) followed by an algorithmic 
solution of the correspondent group-theoretic problem.

Probably, due to preferences to algorithmic solutions, 
it happens that in the projective knot theory the fundamental group has not
been used so far. Of course, some general results are applicable. For
example, Stallings Theorem \cite{St} about necessary and sufficient 
conditions for fibering of a compact 3-manifold over a circle gives a 
criterion for a projective knot to have complement fibered over circle 
in terms of the knot group. But there was no results of this direction 
belonging and specifically to the projective link theory. The conceptual
beauty of fundamental group makes them inevitable.   

In this paper we  
prove three theorems about relation between geometric properties of
a projective link $L$ and properties of $\pi_1(\rpt\sminus L)$. 

\subsection*{Homological classes of the link components}\label{s1.2.5}
The homology group $H_1(\R^3)$ of $\R^3$ is trivial. Therefore each
connected component of a classical link is zero-homologous. 

The homology
group $H_1(\rpt)$ consists of two elements. Therefore we can distinguish
projective knots homologous to zero and non-homologous to zero. In a
projective link, each connected component realizes one of these two classes. 
The distribution of components into two classes has profound effect on other 
invariants and properties of a projective link.
See Lemma \ref{lemHlgy} below. 

\subsection*{Projective unknots}\label{s1.3} 

Recall that a classical knot is called an \text{unknot} 
if it is isotopic to a circle in a plane.

\begin{figure}[htb]
\centerline{
\includegraphics[scale=.65]{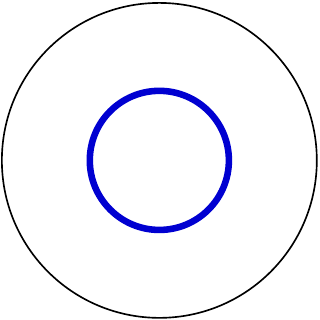}\hspace*{1cm}
\includegraphics[scale=.65]{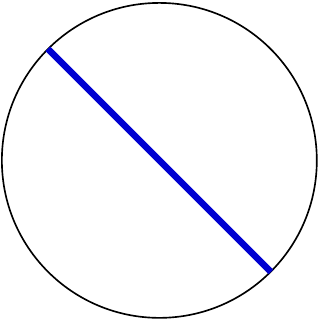}
}
\caption{Diagrams of projective unknots}
\label{f2}
\end{figure}

In the projective space, knots of two isotopy types clearly deserve the 
name of {\em unknot.\/} One of 
the classes consists of projective knots isotopic to unknots in 
$\R^2\subset\rpt$. We call  knots of this kind {\em affine unknots.}
The other isotopy type contains real projective lines. 
Projective diagrams of the unknots are shown  in 
Figure \ref{f2}. Notice that there is one unknot in each of the two
homology classes: affine unknot is zero-homologous, a projective line is
not.    

Dehn's characterization of unknots in $\R^3$
(as classical knots with knot group isomorphic to $\Z$)
 has the following natural counterparts for both types of 
projective unknots. 

\begin{Th}\label{Th1} 
A knot $K\subset\rpt$ is isotopic to the affine unknot if and only if
$\pi_1(\rpt\sminus K)$ is isomorphic to the free product of $\Z$ and
$\ZZ$.
\end{Th}

\begin{Th}\label{Th2} 
A knot $K\subset\rpt$ is isotopic to a projective line if and only if
$\pi_1(\rpt\sminus K)$ is isomorphic to $\Z$. 
\end{Th}

These theorems are proved  in Section \ref{s3}.

\subsection*{Contractibility}\label{s1.4}
A projective link $L$ is said to be {\em contractible\/} if
it is isotopic to a link, which does not intersect the projective plane 
$\rpp\subset \rpt$. 
The complement $\rpt\sminus\rpp$ can be identified 
with the affine space $\R^3$. Therefore a contractible projective link
$L$ can be contracted by a continuous deformation $L_t$,
$t\in[0,1]$, which starts with $L_0=L$ and is an isotopy except for the 
last moment $t=1$, when $L_1$ is a point.
  
In the papers of the first author \cite{D1}, \cite{D2} and the second 
author \cite{V} 
contractible projective links were called {\em affine.\/} 

The following theorem provides necessary and sufficient condition
for a projective link to be contractible:

\begin{Th}\label{Th3}
A link $L\subset \rpt$ is contractible if and only if 
$\pi_1(\rpt\sminus L)$ contains a non-trivial element of order 2. 
\end{Th}

The problem of determining whether a link is contractible was considered in
literature and there are results in this direction, mostly about necessary
conditions:

{\bf Homology condition.}  Each connected component of a link 
$L\subset\rpt$ realizes a homology class, an element of 
$H_1(\rpt;\ZZ)=\ZZ$.   
All components of an contractible link $K\subset\rpt$ realize  
$0\in H_1(\rpt;\ZZ)$.  The converse is not true: there exist knots 
homological to zero in $\rpt$, which are not contractible. 

The homology class of a closed curve in $\rpt$ equals its intersection 
number with a projective plane modulo 2. 
On a projective link diagram, the number of its points on the boundary 
circle is twice the intersection number of the link with a projective
plane. Therefore the
homology class realized by a projective link is non-zero iff the number 
of the boundary points of its diagram is not divisible by four. 

\begin{figure}[ht]
\centerline{\includegraphics[scale=.75]{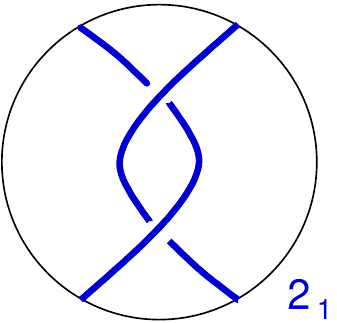}\hspace*{1cm}
\includegraphics[scale=.75]{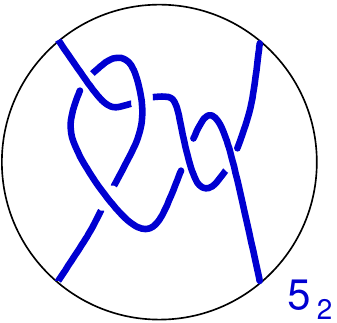}\hspace*{1cm}
\includegraphics[scale=.75]{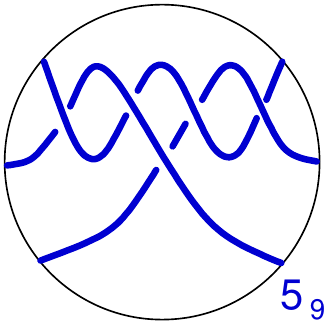}}
\caption{Diagrams of non-contractible links}
\label{f3}
\end{figure}

All three projective knots shown in Figure \ref{f3} are not contractible.
The first two of them, knots $2_1$ and $5_2$ are zero-homologous, because
the numbers of their boundary points equal 4, the third
one, $5_9$, is not zero-homologous, as the number of its boundary points is
6.

{\bf Self-linking number.} If a knot $K\subset\rpt$ realizes 
$0\in H_1(\rpt;\ZZ)$, then a {\em self-linking number\/} $\sl(K)\in\Z$ is 
defined as the linking number modulo 2  of the connected 
components of the preimage $\widetilde K\subset S^3$ of $K$ under the 
covering $S^3\to\rpt$, see \cite{D1}, \S7. 

If $K$ is contractible, then $\sl(K)=0$, see \cite{D1}, \S7. 
For the knot $2_1$ shown in Figure \ref{f3}, this invariant equals 2, 
and this is why $2_1$ is not contractible.

{\bf Exponents of monomials in the bracket polynomial.} 
 If a projective $k$ component link $L$ is contractible, 
then the exponents of all monomials of its bracket polynomial $V_K$ defined 
by the first author in \cite{D1} are congruent to $2k-2$ modulo 4, see 
\cite{D1}, Theorem 7. 

The self-linking and the exponent conditions are independent. 
For the knot $5_2$ the self-linking condition does not work, as $\sl5_2=0$,
while the exponent condition works: $V_{5_2}=A^4+A^2-1-2A^{-2}+A^{-4}+2A^{-6}-2A^{-10}+A^{-14}$. 
On the other hand, for the knot $5_9$ the self-linking condition works, as 
$\sl5_9=3$, while the exponent condition does not: 
$V_{5_9}=A^{-8}+A^{-12}-A^{-20}$. 
See Figure \ref{f3} and the projective link table \cite{D2}. 

In \cite{D1} the first author introduced a notion of alternating projective
link diagram and found the following criterion for contractibility of 
projective links 
which admit alternating diagram: a projective link represented by an
alternating diagram is contractible if and only if the difference between the
greatest and lowest exponents of its bracket polynomial is divisible by 4.

Theorem \ref{Th3}  provides necessary and sufficient condition 
for contractibility of a link in $\rpt$. 
However the criteria formulated in terms of 
fundamental groups are not easy to verify. 
Therefore the known necessary conditions of contractibility mentioned 
above may happen to be more useful for checking if a specific knot is 
contractible. 

\subsection*{Contractibility in general 3-manifold}\label{s1.5}
The property of a projective link of being contractible admits the 
following obvious reformulation:\\
{\it A link $L\subset \rpt$ is contractible if and only if there exists an
embedding $i:D^3\to\rpt$ such that $L\subset i(D^3)$.} 

This property of links in $\rpt$ admits a generalization to
links in an arbitrary 3-manifold. A link $L$ in a 3-manifold $M$ is called 
{\sfit contractible\/} if there exists an embedding $i:D^3\to M$ of the ball
$D^3$ such that $L\subset i(D^3)$.
The following conjecture would provide a generalization of Theorem \ref{Th3} 
to closed 3-manifolds with trivial $\pi_2$.\medskip

\noindent{\bf Conjecture.} {\it Let $M$ be a closed 3-manifold which
satisfies the condition\footnote{As follows from the Papakyriacopoulos sphere theorem, 
this
condition can be reformulated as non-existence of a 2-sphere $\GS$ 
embedded into 
the orientation covering space of $M$ in such a way that $\GS$ does not bound 
a 3-ball there.}
$\pi_2(M)=0$. 
Then a link $L\subset M$ is contractible if and only if $\pi_1(M\sminus L)$
contains a subgroup $G$ which is mapped isomorphically onto 
$\pi_1(M)$ by the inclusion homomorphism $\pi_1(M\sminus L)\to\pi_1(M)$.}

\subsection*{Organization of the paper}\label{s1.7}
In Section \ref{s2} we explain how to write down a presentation of the 
group  $\pi_1(\rpt\sminus L)$ by generators and relations for a projective
link $L$ given by its projective link diagram. 
Our presentation  is similar to the
presentation for the classical knot group introduced by Dehn in
\cite{DH}.  Moreover, our presentation can be considered a 
generalization of slightly enhanced version of the Dehn presentation.
We found beneficial preliminary checkerboard coloring of a link diagram
both in classical and projective environments. 

In the classical knot theory Wirtinger presentations are used more often 
than Dehn presentations. It also generalizes to the projective environment,
but becomes more cumbersome, and we do not consider it.

Section \ref{s3} is devoted to proof of Theorems \ref{Th1} and \ref{Th2}. 
Section \ref{s4} contains a proof of Theorem \ref{Th3}.

\section{Generators and relations for $\pi_1(\rpt\sminus L)$}\label{s2}
The material of this section is not necessary for understanding of subsequent sections.

\subsection*{Geometric origin of generators and relations}\label{s2.1}
The generators of our presentation correspond to the connected components 
of the complement of the link projection. Relations are of two types. 
Relations of the first type correspond to crossings. Relations of the 
second type correspond to pairs of antipodal arcs of the boundary circle. 

For example, the presentation of group for the projective knot $2_1$ in
Figure \ref{f3} has 5 generators, 2 relations of the first type and 2
relations of the second type. 

\subsection*{Projective diagrams}\label{s2.2}
In order to explain how our presentation appears, let us recall that a  
projective link diagram of $L\subset\rpt$ comes from the ball model of 
$\rpt$, that is a presentation of $\rpt$ as a ball $D^3$ with each point on
the boundary sphere identified with its antipode:
$$
\rpt=D^3/_{x\sim -x \ \text{ for any } x\in\p D^3}.
$$
 Choose the corresponding 
mapping $D^3\to\R P^3$ so that the image
of the northern and southern poles $N$ and $S$ of the ball does not 
belong to  $L$. Denote by $L'$ 
the preimage of $L$ in $D^3$. Let $p: L'\to D^2$ be the projection 
onto the equatorial disk $D^2\subset D^3$ given by the formula $p:x\mapsto
c(x) \cap D^2$, where $c(x)$ is
the (metric) circle in $\R^3\supset D^3$ passing through 
$x \in L'$ and the poles $N$ and $S$  of  $D^3$. 
\begin{figure}[thb]
\centerline{\begin{tikzpicture}[scale=.8,>=latex]   
\draw[thick](0,0)circle(4);
\draw[thick](0,0)ellipse (4cm and .8cm);
\filldraw(0,3.82)circle(2pt);
\node at (0,4.4){$N$};
\filldraw(0,-3.82)circle(2pt);
\node at (0,-4.4){$S$};
\draw(0,-3.82) arc(-38:38:6.2cm);
\draw[very thick,blue](.2,2.3) arc(-120:-70:2cm);
\draw[very thick,blue](.36,-.2) arc(.-100:-60:3.2cm);
\node at (1.8,2.5){$L'$};
\node at (2.85,-.15){$p(L')$};
\node[rotate=-50] at (-3,0){$D^2$};
\filldraw(.96,2.05)circle(2pt);
\filldraw(1.32,-.22)circle(2pt);
\node at (.75,1.75){$x$};
\node at (.78,0.1){$p(x)$};
\node at (1.45,-2.01){$c(x)$};
\draw[very thick,->](1.30,.5)--(1.30,.3);
\draw[very thick,->](1.135,-1.5)--(1.17,-1.3);
\end{tikzpicture}
}
\caption{Projection to $D^2$ in a ball model of \b{$\rpt$}}
\label{proj2disk}
\end{figure}
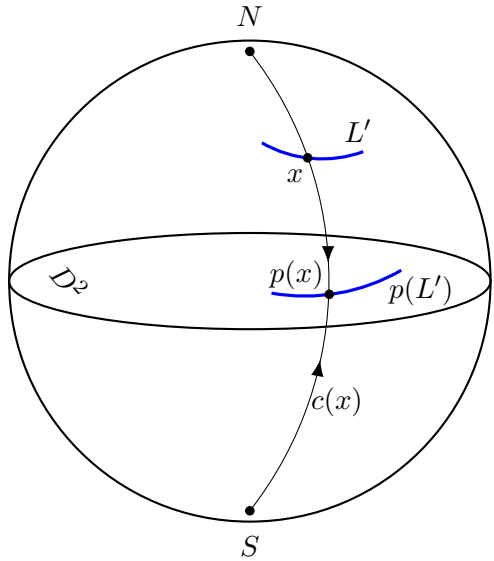

Then we convert $p(L')$ in a diagram by decorating it as usual near each 
double point. 

\subsection*{Choice of a base point}\label{s2.3}
In order to speak about the fundamental group of $\rpt\sminus L$ in detail, 
we have to choose a base point in $\rpt\sminus L$.
Denote by $P$ the point in $\rpt$ obtained from the poles $N$ and $S$ of 
$D^3$. In other words, let $P=p(N)=p(S)\in\rpt$.  


The group $\pi_1(\rpt\sminus L,P)$ consists of homotopy classes of loops 
starting in the base point. Observe, that any path in $D^3$ connecting 
$N$ and $S$ gives rise to a loop like this and hence to an element of 
$\pi_1(\rpt\sminus L,P)$.

\subsection*{Checkerboard coloring and generators}\label{s2.4}
By {\em domains\/} of the diagram, we mean connected components of
$D^2\sminus p(L')$. Choose a checkerboard coloring of domains. 

In each domain choose a point and draw a circular arc connecting
the poles $N$ and $S$ through this point. 
For dark domains orient the arc from $N$ to $S$, for light domains orient 
the arc from $S$ to $N$. This gives 
for each domain a path. For a dark domain the path in $D^3$ starts at $N$
and finishes at $S$, for a light domain the path starts at $S$ and finishes
at $N$. In $\rpt$ the paths gives rise to loops at the base point $P$ and,
further, to elements of $\pi_1(\rpt\sminus L,P)$.

Assign to each domain on the diagram a symbol, and denote the corresponding 
element of $\pi_1(\rpt\sminus L,P)$ by this symbol. 

\begin{lem}\label{lemGen} 
The elements of   $\pi_1(\rpt\sminus L,P)$ corresponding to domains
generate   $\pi_1(\rpt\sminus L,P)$. 
\end{lem}

We leave the proof of Lemma \ref{lemGen} to the reader. It is a standard
exercise, which can be done in many ways using your favorite technique for
calculation of a fundamental group. For example, you can use 
Seifert - van Kampen Theorem, or transversality theorem.

\subsection*{Relations corresponding to crossings}\label{s2.5}
\begin{figure}[h]
\centerline{\includegraphics{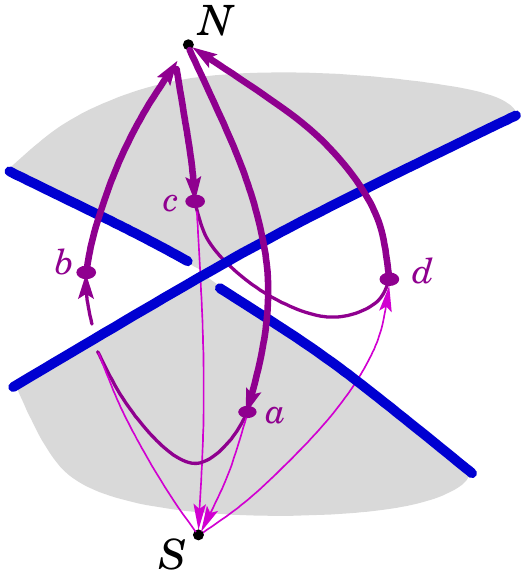}}
\caption{Relation $abcd=1$ at a crossing}
\label{f4.5}
\end{figure}

A relation corresponding to a crossing starts with the symbol of an adjacent 
dark domain, then proceeds with the symbol of the domain adjacent to this 
first domain and separated from it by an overcrossing arc. Then proceeds
with the
symbols of the remaining dark domain followed by remaining light domain. 
All four domains are adjacent to the crossing and are going all the way 
around the crossing. See Figure \ref{f4.5}. 

\begin{figure}[h]
\centerline{\includegraphics[scale=.75]{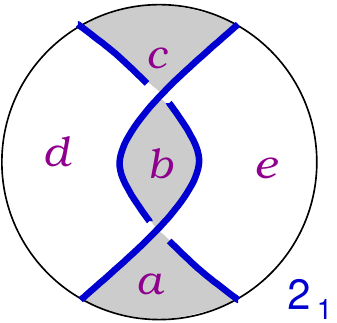}}
\caption{Relations at crossings of $2_1$}
\label{f5}
\end{figure}

On Figure \ref{f5} relations at crossings of a diagram for $2_1$ are shown.
We assign relation $adbe=1$ to the lower crossing if we start with 
the domain $a$, or $bead=1$, if we start at the domain $b$. 
It does not matter where to start, as the relations are
equivalent. Similarly, at the upper crossing we get relation $bdce=1$.


\subsection*{Relations corresponding to pairs of antipodal boundary 
arcs}\label{s2.6}
If the domains adjacent to a pair of antipodal boundary arcs have {\em the 
same color\/}, then relation correspondent to the pair of arcs says that 
the generators of the adjacent domains are {\em inverse\/} to each other. 

On Figure \ref{f5}, all the pairs of antipodal boundary arcs are adjacent to 
domains of  the same color. Thus, we have relations $c=a^{-1}$ and
$e=d^{-1}$. The whole presentation of the group looks as follows:
\begin{multline*}
\{a,b,c,d,e\mid adbe=1, bdce=1, c=a^{-1}, e=d^{-1}\}\\=
\{a,b,d\mid adbd^{-1}=1, bda^{-1}d^{-1}=1\}\\
=\{a,b,d\mid a=db^{-1}d^{-1}=d^{-1}bd\}
=\{b,d\mid d^2=bd^2b\}
\end{multline*}

If the domains adjacent to the pair of antipodal boundary arcs have
{\em different colors\/}, then relation correspondent to the pair of arcs says that 
the generators of the adjacent domains are {\em equal\/} to each other. 
In this case, it's convenient to take into account the relations of this
kind on the preceding step, i.e., while assigning symbols of generators to
the domains. Of course, this can be done in the other case, too. Given a 
projective link diagram equipped with checkerboard coloring and symbols, 
one can write down the correspondent presentation of the group by listing 
the symbols and writing down the 4-term relations at each crossing. 
\begin{figure}[h]
\centerline{
\includegraphics[scale=.75]{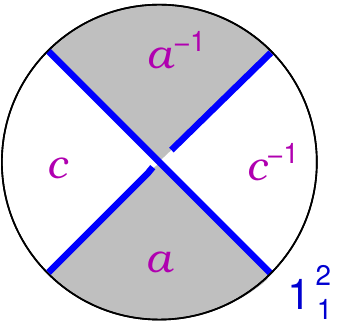}\hspace*{1cm}
\includegraphics[scale=.75]{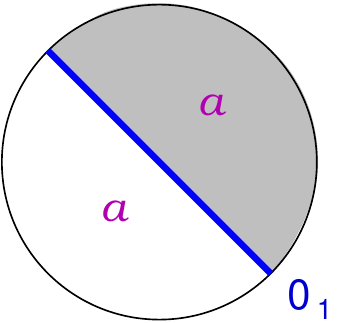}\hspace*{1cm}
\includegraphics[scale=.75]{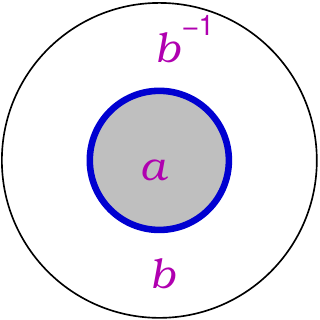}
}
\caption{The graphic stage of writing down the group representation 
}
\label{f6}
\end{figure} 

The simplest complete examples of writing down the group presentations are
shown in Figure \ref{f6}. 
The group for the projective link $1_1^2$, which is a pair of two skew
lines, is $\{a,c\mid ac^{-1}a^{-1}c=1\}$, that is a free abelian group of
rank two. 

The group for $0_1$, a projective line, has generator $a$ and no
relations, so this is an infinite cyclic
group. 

The group of the last knot is $\{a,b\mid
b=b^{-1}\}=\{a,b\mid b^2=1\}$, that is the free product $\Z*\ZZ$ of an 
infinite cyclic group by a group of order two.

\subsection*{A version of Dehn presentation for classical link}\label{s2.4}
If a projective link $L$ is presented by a diagram which does not meet 
the boundary circle of its disk, the $L$ does not intersect the projective
plane which is the quotient space of the boundary sphere $\p D^3$. 
Then removing the boundary circle turns the projective diagram into a 
diagram of the corresponding classical link. 

A presentation of the group $\pi_1(\R^3\sminus L)$ can be obtained from the
presentation described above by removing the generator which corresponds to 
the external domain of the initial projective diagram. This generator had
order two in the group of projective link. Cf. Theorem \ref{Th3}.

\section{Projective unknots}\label{s3}
 
In this section we prove Theorems \ref{Th1} and \ref{Th2}. Recall their
statements.\vspace{8pt}

\noindent{\bf Theorem 1.} 
{\it A knot $K\subset\rpt$ is isotopic to the affine unknot if and only if
$\pi_1(\rpt\sminus K)$ is isomorphic to the free product of $\Z$ and
$\ZZ$.}\vspace{8pt}

\noindent{\bf Theorem 2.} 
{\it A knot $K\subset\rpt$ is isotopic to a projective line if and only if
$\pi_1(\rpt\sminus K)$ is isomorphic to $\Z$.}\vspace{8pt}

 Both theorems admit a proof similar to the classical proof of Dehn's Theorem. 

The "only if" parts are proved on Figure \ref{f6}, where
the groups are found for the standard diagrams of the unknots. 
The rest of this section is devoted to proof of the "if" part.

\subsection*{Plan of the proof}\label{s3.0}
Let $N$ be a tubular neighborhood of $K$, and $T$ be its boundary $\p N$. 
Since $\rpt$ is orientable, $N$ is homeomorphic to $S^1\times D^2$, and 
$T$ is torus, $T$ is homeomorphic to $S^1\times S^1$.

First, we prove that the inclusion homomorphism 
$\pi_1(T)\to\pi_1(\rpt\sminus K)$ 
cannot be injective in either cases. 

Second, by a version of Papakyriacopoulos' Loop Theorem, 
we conclude that there
exists an embedded circle in $T$ which does not bound in $T$, but
bounds a smoothly embedded disk $\GD$ in $\rpt\sminus\Int N$.

Third, we find the homology class realized by the boundary of $\GD$
in $H_1(T)$. The inclusion homomorphism $H_1(T)\to H_1(N)=H_1(K)$
maps this class to a generator in the case of Theorem \ref{Th1} and to the
double of a generator in the case of Theorem \ref{Th2}.

Fourth, the disk $\GD$ is completed by an annulus in $N$ into a disk 
bounded by
$K$, in Theorem \ref{Th1}, and into a projective plane in which $K$ is
contained and does not bound, in Theorem \ref{Th2}.

\subsection*{Non-contractible loop on $T$ contractible in
$\rpt\sminus K$}\label{s3.1}
The assumptions of each of 
the theorems  imply that  the inclusion homomorphism 
$\pi_1(T)\to\pi_1(\rpt\sminus K)$ is not injective.

For Theorem \ref{Th2} it is obvious: there is no monomorphism 
$\Z\times \Z\to\Z$.

For Theorem \ref{Th1} it is also true, but requires more arguments.

\begin{lem}\label{lemGr}
There is no monomorphism $\Z\times\Z\to\Z*\ZZ$. 
\end{lem}
\begin{proof}Due to Kurosh Subgroup Theorem \cite{Kur}, 
a subgroup of $\Z*\ZZ$ should be a free product of cyclic groups. 
A group isomorphic to $\Z\times\Z$ cannot
be a free product of cyclic groups by Baer-Levi Theorem \cite{BL}.
\end{proof}

\subsection*{Disk provided by the Loop Theorem}\label{s3.2}
Recall that according to the Loop Theorem (see, e.g., \cite{Hatcher}), if 
$M$ is a 3-manifold, and the inclusion homomorphism $\pi_1(\p M)\to\pi_1(M)$
is not injective, then there exists an embedding $f:(D^2,\p D^2)\to(M,\p M)$
such that it gives a non-contractible loop $\p D^2\to\p M$.

\begin{lem}\label{lemLoopTh}
Under assumptions of Theorem \ref{Th1} or of Theorem \ref{Th2}, let $N$ be
a tubular neighborhood of $K$ and let $E$ be the exterior $\rpt\sminus\Int
N$ of $K$. Then there exists a disk \b{$\GD$} properly embedded in 
$\rpt\sminus\Int N$ such that its boundary $\GS=\p\GD$ is not contractible
in $\p N$.
\end{lem}

\begin{proof}
Apply the Loop Theorem to $M=E$. It gives an embedding 
$f:(D^2,\p D^2)\to(E,\p N)$ such that its restriction 
$\p D^2\to \p N$ is a non-contractible loop. 
Let us denote the disk $f(D^2)$ by $\GD$ and the simple closed curve 
$f(\p D^2)$ by $\GS$.
\end{proof}

A simple closed curve on a torus non-homological to zero realizes an
indivisible homology class. Thus, $\GS$ realizes an indivisible class
in $H_1(T)$. 

\subsection*{Homological calculation}\label{s3.3}
There are exactly two elements in $H_1(\rpt)=\ZZ$: the zero and non-zero. 

\begin{lem}\label{lemHlgy}
Let  $K$ be a projective knot. Then 
$$
H_1(\rpt\sminus K)=
\begin{cases} \Z\oplus \ZZ &\text{if $K$ is zero-homologous in $\rpt$}\\
 \Z &\text{if $K$ is homologous to a projective line}\end{cases}
$$
$$
H_1(\rpt,K)=
\begin{cases} \ZZ &\text{if $K$ is zero-homologous in $\rpt$}\\
0   &\text{if $K$ is homologous to a projective line}
\end{cases}
$$
The group $H_2(\rpt,K)=\Z$ independently on the homology class realized by
$K$, but the boundary homomorphism $H_2(\rpt,K)\to H_1(K)$ depends: if the
homology class of $K$ is zero, then this is an isomorphism, if not, then 
this is a monomorphism with the image of index two. 
\end{lem}

\begin{proof}
Consider the homology sequence of pair $(\rpt,K)$:
$$
\begin{tikzcd}[cramped,sep=small] 
H_2(K)\arrow[r]\arrow[d,equal]
&H_2(\rpt,K)\arrow[r]\arrow[d,equal]
&H_1(K)\arrow[r]\arrow[d,equal]
&H_1(\rpt)\arrow[r]\arrow[d,equal]
&H_1(\rpt,K)\arrow[r]\arrow[d,equal]
&\tilde H_0(K)\arrow[d,equal]\\
0\arrow[r]
&?\arrow[r]
&\Z\arrow[r]
&\ZZ\arrow[r]
&?\arrow[r]
&0
\end{tikzcd}
$$
If the middle inclusion homomorphism $H_1(K)\to H_1(\rpt)$ is zero, then
by exactness of the sequence we get $H_1(\rpt,K)=\ZZ$ and  $H_2(\rpt,K)=\Z$ .
If the homomorphism is not zero, then $H_1(\rpt,K)=0$ and still 
$H_2(\rpt,K)=\Z$, but the boundary homomorphism $H_2(\rpt,K)\to H_1(K)$ 
  in the first case is an isomorphism and in the second case it is a
monomorphism with the image $2\Z\subset \Z$.

Further, by duality and the universal coefficient formula
$$H_1(\rpt\sminus
K)=H^2(\rpt,K;\Z)=\Hom(H_2(\rpt,K),\Z)\oplus\Ext(H_1(\rpt,K),\Z)$$ 
Hence, 
$$
H_1(\rpt\sminus K)=
\begin{cases} \Z\oplus \ZZ &\text{if $K$ is zero-homologous in $\rpt$}\\
 \Z &\text{if $K$ is homologous to a projective line}\end{cases}
$$
\end{proof}
Since $K$ is either zero-homologous or homologous to a projective line,
we see that the group $H_1(\rpt\sminus K)$ determines which homology class 
is realized by $K$. Hence, under assumptions of Theorem \ref{Th1}, $K$ is
zero homologous, while, under assumptions of Theorem \ref{Th2}, $K$ is
homologous to a projective line.

\subsection*{Mutual position of $K$ and $\GS$ in $N$}\label{s3.4} 
\begin{lem}\label{lemInN} Under the assumptions of Theorem \ref{Th1}, the curves $K$
and $\GS$ bound an annulus embedded in $N$. Under assumptions of Theorem
\ref{Th2}, the curve $\GS$ is a boundary of a M\:obius band embedded in
$N$.
\end{lem}

\begin{proof}
Recall that by Lemma \ref{lemLoopTh}, under assumptions of Theorem \ref{Th1} or
Theorem \ref{Th2}, there is a simple closed curve $\GS$ on the 
boundary $T$ of a tubular neighborhood $N$ of $K$.

Denote by $E$ the exterior $\rpt\sminus\Int N$ of the knot $K$.
 The curve $\GS$ bounds a disk
$\GD$ in $E$ and is not zero-homologous on $T$. Denote by $\Gs\in H_1(T)$ 
the homology class realized by $\GS$ and by $\Gd\in H_2(E,T)$ the homology class 
realized by the disk $\GD$. 

The knot $K$ is a deformation retract of its tubular neighborhood $N$.
Hence the pair $(\rpt,K)$ is homotopy equivalent to $(\rpt,N)$,
the inclusion homomorphism $H_2(\rpt,K)\to H_2(\rpt,N)$ is an isomorphism
and, by Lemma \ref{lemHlgy}, $H_2(\rpt,N)=H_2(\rpt,K)=\Z$. 
By excision, $H_2(E,T)=H_2(\rpt,N)=\Z$. So, we have the following
commutative diagram:
$$
\begin{tikzcd} 
\Z\arrow[r,equal]&H_2(\rpt,K)\arrow[d,swap,"\text{iso}"]\arrow[r]&H_1(K)\arrow[d,"\ink_*
\text{ iso}"]\arrow[r,equal]
&\Z\\
\Z\arrow[r,equal]&H_2(\rpt,N)\arrow[r]&H_1(N)\arrow[r,equal]
&\Z\\
\Z\arrow[r,equal]&H_2(E,T)\arrow[u,"\text{iso}"]\arrow[r]&H_1(T)\arrow[u,swap,"\ink_*"]\arrow[r,equal]&\Z\times\Z
\end{tikzcd}
$$

Since $\Gs\ne0$ and $GS$ is a simple closed
curve, $\Gs$ is indivisible in $H_1(T)=\Z\times\Z$. Hence $\Gd$ is
indivisible in $H_2(E,T)=\Z$, so $\Gd$ generates $H_2(E,T)$. The vertical
homomorphism on the left hand side of the diagram map $\Gs$ to a generator
of $H_2(\rpt,K)$. By Lemma \ref{lemHlgy}, this generator is mapped by 
the top homomorphism $H_2(\rpt,K)\to H_1(K)$ under assumptions of Theorem
\ref{Th1} to a generator of $H_1(K)$, or, under assumptions of Theorem
\ref{Th2} to a doubled generator of $H_1(K)$.
\end{proof}

\subsection*{Completion of the proofs}\label{s3.5}
By Lemma \ref{lemInN}, under assumptions of Theorem \ref{Th1}, there exists 
an annulus $A\subset N$ bounded by $\GS\cup K$. The union $A\cup \GD$ is 
homeomorphic to disk. It is bounded by $K$. By an isotopy fixed on $K$ it 
can be made smoothly embedded in $\rpt$. One can contract $K$ along this 
disk to a small planar circle. This completes our proof of Theorem
\ref{Th1}.\qed

By Lemma \ref{lemInN}, under assumptions of Theorem \ref{Th2},
$\GS$ is a boundary of a M\:obius band $M$ embedded in $N$ and $K$ lies on
$M$ as the midline. The union $\GD\cup M$ is homeomorphic to $\rpp$. 
The original knot $K$ lies on it as a projective line. By an isotopy the union 
can be made smooth in $\rpp$. Any smoothly embedded $\rpp$ in
$\rpt$ is smoothly isotopic to a projective embedding. 
Then $K$ can be deformed by a smooth isotopy to
a projective line. This completes our proof of Theorem \ref{Th2}.\qed

\section{Contractibility}\label{s4} 
\subsection*{Proof of Theorem \ref{Th3}. Necessity.}\label{s4.1} 
Assume that $L$ is an contractible link. Since 
the group $\pi_1(\rpt\sminus L)$ is invariant under isotopy, it suffices
to prove that $\pi_1(\rpt\sminus L)$ contains an element of order 2 if 
$L$ does not intersect a plane $\rpp\subset\rpt$.

If $L\cap\rpp=\empt$, then $L$ is contained in the contractible part 
$\rpt\sminus\rpp=\R^3$ of $\rpt$,
and  $\rpt\sminus L$ is a union of $\rpt\sminus(L\cup\rpp)$ and a regular 
neighborhood of $\rpp$. The van Kampen Theorem applied to this presentation 
of $\rpt\sminus L$ implies that  $\pi_1(\rpt\sminus L)$ is the free product
$\ZZ*\pi_1(\R^3\sminus L)$. Hence it contains a  non-trivial element of
order 2.\hfill\qed\bigskip

 \subsection*{A lemma about a map of the projective plane}\label{s4.2}
The proof of sufficiency is prefaced with the following simple 
homotopy-theoretic lemma: 

\begin{lem}\label{lem1} Let $f:\rpp\to\rpp$ be a map inducing isomorphism
on fundamental group. Then the covering map $\widetilde f:S^2\to S^2$ is
not null homotopic.
\end{lem}

\begin{proof} Consider the diagram
$$
\minCDarrowwidth8pt
\begin{CD}
0@>>>H_2(\rpp;\ZZ)@>>>H_2(S^2;\ZZ)@>>>H_2(\rpp;\ZZ)@>>>H_1(\rpp;\ZZ)@>>>0\\
@. @VV{f_*}V @VV{\widetilde f_*}V @VV{f_*}V  @VV{f_*}V @. \\
0@>>>H_2(\rpp;\ZZ)@>>>H_2(S^2;\ZZ)@>>>H_2(\rpp;\ZZ)@>>>H_1(\rpp;\ZZ)@>>>0
\end{CD}
$$ 
in which the rows are segments of the Smith sequence (see, e.g., \cite{Bre}) 
for the antipodal involution on $s:S^2\to S^2:x\mapsto-x$. The diagram 
is commutative, because $\widetilde f$ commutes with $s$. 
Notice that all the groups in this diagram are isomorphic to $\ZZ$. 
Exactness of the Smith sequences implies that the middle horizontal 
arrows in both rows are trivial, and the horizontal arrows next to them
are isomorphism. By assumption, the rightmost vertical arrow is an
isomorphism. Therefore the next vertical arrow is an isomorphism. This
isomorphism 
coincides with the homomorphism represented by leftmost vertical arrow.
Hence, the next arrow $\widetilde f_*:H_2(S^2;\ZZ)\to H_2(S^2;\ZZ)$ is 
an isomorphism. Thus $\widetilde f$ is not null homotopic.    
\end{proof}

\subsection*{Proof of Theorem \ref{Th3}. Sufficiency}\label{s4.3}
Assume that $\pi_1(\rpt\sminus L)$ contains a
non-trivial element $\Gl$ of order two. Realize $\Gl$ by a smoothly 
embedded loop
$l:S^1\to \rpt\sminus L$. 
\subsubsection*{From a loop of order 2 to a singular projective
plane in $\rpt\sminus L$}\label{s4.3.1}
Since $\Gl^2=1$, there exists a
continuous map $D^2\to \rpt\sminus L$ such that its restriction to the
boundary circle $\p D^2$ is the square of $l$. Together with $l$, this 
map gives a continuous map of the projective plane 
$P=D^2/_{x\sim -x, \text{ for } x\in\p D^2}$ to $\rpt\sminus L$. 
%
Denote this map by $g$. 
So, $g:P\to\rpt\sminus L$ is a generic differentiable map 
which induces
a monomorphism $g_*$ of $\pi_1(P)=\ZZ$ to $\pi_1(\rpt\sminus L)$. 

\subsubsection*{Singular sphere over the singular projective plane}\label{s4.3.2}  
Let $p:S^3\to\rpt$ be the canonical two-fold covering. Consider its
restriction $S^3\sminus p^{-1}(L)\to\rpt\sminus L$. 
Observe that $\Gl$ does not belong to the group of this covering, because 
otherwise $\pi_1(S^3\sminus p^{-1}(L))$ would contain a non-trivial element 
of order two, which is impossible - a link group does not have any
non-trivial element of finite order. 

Therefore the covering of the projective plane induced from $p$ via $g$
is a non-trivial two-fold covering. Its total space is a 2-sphere. Denote
it by $S$. The map $\widetilde g$ which covers $g$ maps $S$ to    
$S^3\sminus p^{-1}(L)$. 

\subsubsection*{Non-contractibility of the
singular sphere in $S^3\sminus p^{-1}(L)$}\label{s4.3.3}
Let us choose a point $x\in L$.
Its complement $\rpt\sminus\{x\}$ is homotopy equivalent to $\rpp$. 
Indeed, the projection from $x$ to any projective plane, which does not 
contain $x$, is a deformation retraction  $\rpt\sminus\{x\}\to\rpp$. 
The composition of
$g:P\to\rpt\sminus L$ with the inclusion $\rpt\sminus L\to \rpt\sminus\{x\}$
induces an isomorphism $\pi_1(P)\to\pi_1(\rpt\sminus\{x\})$. Both spaces, 
$P$ and $\rpt\sminus\{x\}$, have homotopy type of $\rpp$. 
Lemma \ref{lem1} implies that $\widetilde g:S\to S^3\sminus p^{-1}(x)$ 
is not null-homotopic. 

\subsubsection*{Existence of a non-singular sphere $\GS_0\subset S^3$
that splits $p^{-1}(L)$}\label{s4.3.4}
Denote by $\Gs$ the antipodal involution $S^3\to S^3:x\mapsto -x$.

\begin{lem}\label{lem2} There exists a non-singular polyhedral 
submanifold $\GS_0$ of $S^3$
homeomorphic to $S^2$ such that $\GS_0\cap p^{-1}(L)=\empt$ and the two points
of $p^{-1}(x)$ belong to different connected components of
$S^3\sminus\GS_0$.
\end{lem}

\begin{proof} First, let us apply Whitehead's modification \cite{Whi} of 
the Papakyriakopoulos Sphere Theorem \cite{Papa}. 

Recall the statement of this theorem (Theorem (1.1) of \cite{Whi}): 
{\it For any connected,
orientable triangulated 3-manifold $M$ and subgroup $\GL\subset\pi_2(M)$  
which is invariant under the action of $\pi_1(M)$, if $\GL\ne\pi_2(M)$, 
then $M$ contains a non-singular polyhedral 2-sphere which is essential 
$\mod\GL$.} 

We will apply this theorem to $$M=S^3\sminus p^{-1}(L), \; \; 
\GL=\Ker\left(\ink_*:\pi_2(S^3\sminus p^{-1}(L))\to\pi_2(S^3\sminus
p^{-1}(x))\right).$$
We know that $\pi_2(S^3\sminus\{x\})=\pi_2(S^2)=\Z$ and that 
the homotopy class of $\widetilde g$ is non-trivial in 
$\pi_2(S^3\sminus p^{-1}(x))$. Therefore the homotopy class of $\widetilde
g$ does
not belong to $\GL$, and hence $\GL\ne\pi_2(S^3\sminus p^{-1}(L))$. Thus, 
all the assumptions of the Whitehead theorem are fulfilled.  

Let us denote by $\GS_0$  a non-singular polyhedral 2-sphere whose
existence is stated by the Whitehead theorem. By the Alexander Theorem
\cite{Alx}, $\GS_0$ bounds in $S^3$ two domains homeomorphic to ball. 
Since $\GS_0$ is not null homotopic in $S^3\sminus
p^{-1}(x)$, each of these domains contains a point of $p^{-1}(x)$.
\end{proof}

\subsubsection*{Improving the splitting sphere }\label{s4.3.5}
\begin{lem}\label{lem3} There exists a smooth submanifold $\GS$ of $S^3$
homeomorphic to $S^2$ such that $\GS\cap p^{-1}(L)=\empt$, the two points
of $p^{-1}(x)$ belong to different connected components of
$S^3\sminus\GS$,
and either $\GS=\Gs(\GS)$ or $\GS\cap\Gs(\GS)=\empt$. 
\end{lem}

\begin{proof}
Any polyhedral compact surface can be approximated by a smooth
2-submanifold. Let  $\GS_1$ be a smooth submanifold of $S^3$ approximating
$\GS_0$ in $S^3\sminus p^{-1}(L)$. 

The antipodal involution $\Gs:S^3\to S^3$ is an
automorphism of the covering $p:S^3\to\rpt$, therefore $p^{-1}(L)$ is
invariant under $\Gs$ and $\Gs(\GS_1)\subset S^3\sminus p^{-1}(L)$.

Let us assume that $\GS_1$ and $\Gs(\GS_1)$ are transversal -- this can be
achieved by an arbitrarily small isotopy of $\GS_1$. 
Then the intersection
$\GS_1\cap\Gs(\GS_1)$ consists of disjoint circles. 


Take a connected component $C$ of $\GS_1\cap\Gs(\GS_1)$ which is innermost in
$\Gs(\GS_1)$ 
(i.e., which  bounds in $\Gs(\GS_1)$ a disc $D$ containing no other 
components of $\GS_1\cap\Gs(\GS_1)$). 

First, assume that $C\ne \Gs(C)$. In this case, make surgery on $\GS_1$ 
along $D$:
remove a regular neighborhood $N$ of $C$ from $\GS_1$ and attach to $\p N$ 
two discs parallel to $D$. This surgery does not change the homology class 
with coefficients in $\ZZ$ realized by $\GS_1$ in $S^3\sminus p^{-1}(x)$.

Denote by $\GS_2$ the result of this surgery on $\GS_1$. This is a disjoint 
union of two spheres. 
The sum of the homology classes realized by them is the same 
non-trivial element
of $H_2(S^2\sminus p^{-1}(x);\ZZ)=\ZZ$ which was realized by $\GS_1$. 
Therefore, one of the summands is
non-trivial. The corresponding component of $\GS_2$ separates the two points 
of $p^{-1}(x)$. Denote this component by $\GS_3$.
 Since $C\ne\Gs(C)$, the number of connected components of
$\GS_3\cap\Gs(\GS_3)$ is less than the number of connected components of 
$\GS_1\cap\Gs(\GS_1)$, all other properties of $\GS_1$ are inherited by
$\GS_3$, and we are ready to continue with the next connected component
of $\GS_3\cap\Gs(\GS_3)$ which bounds in $\Gs(\GS_3)$ a disc containing 
no other components of $\GS_3\cap\Gs(\GS_3)$).

Second, consider the case $C=\Gs(C)$.  Then the disc $D\subset\Gs(\GS_1)$
together with it's image $\Gs(D)\subset \GS_1$ form an embedded sphere, which
is invariant under $\Gs$ and does not meet the rest of $\GS\cup\Gs(\GS_1)$
besides along $C$, that is
$(D\cup\Gs(D))\cap((\GS_1\sminus\Gs(D)\cup(\Gs(\GS_1)\sminus D))=\empt $. 

If $D\cup\Gs(D)$ separates points of $p^{-1}(x)$, we are done: 
we can smoothen the corner of  $D\cup\Gs(D)$ along $C$ keeping it invariant
under $\Gs$ and take the result for $\GS$. 

If $D\cup\Gs(D)$ does not separate points of $p^{-1}(x)$, then 
$D\cup(\GS_1\sminus\Gs(D))$  separates points of $p^{-1}(x)$ 
(as well as its image under $\Gs$, that is 
$\Gs(D)\cup(\Gs(\GS)\sminus D)$). Indeed, the homology classes realized by 
$D\cup\Gs(D)$ and $D\cup(\GS_1\sminus\Gs(D))$ in $S^3\sminus p^{-1}(x)$ 
differ from each other by the homology class of
$\Gs(D)\cup(\GS_1\sminus\Gs(D))=\GS_1$ which is known to be nontrivial. 
So, if the class of  $D\cup\Gs(D)$ is trivial, then the class of 
$\Gs(D)\cup(\Gs(\GS_1)\sminus D)$ is not. 
Then smoothing of a corner along $C$ turns $D\cup(\GS_1\sminus\Gs(D))$ 
into a new sphere $\GS_2$ such that
$\GS_2\cap\Gs(\GS_2)$ has less connected components than  $\GS_1\cap\Gs(\GS_1)$. 

By repeating this construction, we will eventually build up a sphere 
$\GS\subset S^3\sminus p^{-1}(L)$ with the required properties.
\end{proof}

\subsubsection*{Completion of the proof}\label{s4.3.6}
Let us return to the proof of Main Theorem. If the sphere $\GS$
provided by Lemma \ref{lem3} is invariant under $\Gs$, then $\GS$ divides
$S^3$ into two balls which are mapped by $\Gs$ homeomorphically to each
other. Let $B$ be one of them. The part of $p^{-1}(L)$ contained in 
$B$ can be moved by an isotopy fixed on a neighborhood of the boundary 
of $B$ inside an arbitrarily small metric ball in $S^3$. Using $\Gs$,
extend this isotopy to a $\Gs$-equivariant isotopy of the whole $S^3$.
The equivariant isotopy defines an isotopy of $\rpt$ which moves $L$
to a link contained in a small metric ball. This proves that $L$ is
an contractible link.

Consider now the case in which the sphere $\GS$ provided by Lemma 
\ref{lem3} is not $\Gs$-invariant, but rather is disjoint from its image
$\Gs(\GS)$. Then spheres $\GS$ and $\Gs(\GS)$ divide $S^3$ into three
domains: two of them are balls bounded by $\GS$ and $\Gs(\GS)$,
respectively. Let us denote by $B$ the ball bounded by $\GS$, then
its image $\Gs(B)$ is bounded by $\Gs(\GS)$. Denote the third domain
by $E$. It is invariant under $\Gs$. 

If one of the points from $p^{-1}(x)$ belonged to $E$, then the other one
also would belong to $E$, and then the sphere $\GS$ would be contractible
in $S^3\sminus p^{-1}(x)$. Therefore $B\cap p^{-1}(L)\ne\empt$. Denote 
$B\cap p^{-1}(L)$ by $K$. This is a sublink of $p^{-1}(L)$.
It  can be moved by an isotopy fixed on a neighborhood of the boundary 
of $B$ inside an arbitrarily small metric ball in $S^3$. Then this isotopy 
can be extended to $\Gs$-equivariant isotopy of $S^3$ fixed on $E$.
This equivariant isotopy defines an isotopy of $\rpt$ which moves $p(K)$
to a link contained in a small metric ball. 

Thus our link $L$ is presented as a disjoint sum of an contractible link $p(K)$
and the rest $L\sminus p(K)$ of $L$. 
If $L\sminus p(K)=\empt$, then we are done. If not,
then $\pi_1(S^3\sminus L)$ is presented as a free product of 
$\pi_1(B\sminus K)$ and $\pi_1(\rpt\sminus (L\sminus p(K))$. By the
assumption, the group $\pi_1(S^3\sminus L)$ has a non-trivial element of
order 2. The first factor, $\pi_1(B\sminus K)$ cannot contain such an
element, because this is a group of a classical link. Hence, the second 
factor,  $\pi_1(\rpt\sminus (L\sminus p(K))$, contains it, and we can
apply the constructions and arguments above to the link $L\sminus p(K)$. 
This link contains less
components than the original one, therefore, after several iterations, 
we will come to the 
situation in which $p^{-1}(L)\cap B=\empt$. 
\qed

\end{document}